\documentclass[12pt]{amsart}

\usepackage{amssymb}	
\usepackage{epsfig}		
\usepackage{amsthm}	

\newcommand{\lcr}{\raisebox{-5pt}{\mbox{}\hspace{1pt}
                  \includegraphics{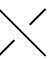}\hspace{1pt}\mbox{}}}
\newcommand{\ift}{\raisebox{-5pt}{\mbox{}\hspace{1pt}
                  \includegraphics{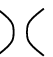}\hspace{1pt}\mbox{}}}
\newcommand{\zer}{\raisebox{-5pt}{\mbox{}\hspace{1pt}
                  \includegraphics{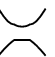}\hspace{1pt}\mbox{}}}

\newcommand{\ot}{\otimes}

\newtheorem{theorem}{Theorem}
\newtheorem{lemma}{Lemma}

\title{Traces on the Skein Algebra of the Torus}

\author{Michael McLendon}
\address{Department of Mathematics and Computer Science,
Washington College, Chestertown, Maryland, 21620}
\email{mmclendon2@washcoll.edu}

\begin{document}

\begin{abstract}
For a surface $F$, the Kauffman bracket skein module
of $F \times [0,1]$, denoted $K(F)$,
admits a natural multiplication which makes it an algebra.  When
specialized at a complex number $t$, nonzero and not a root of unity,
we have $K_t(F)$, a vector space over $\mathbb{C}$.
In this paper, we will use the product-to-sum formula of Frohman and Gelca
to show that the vector space $K_t(T^2)$ has
five distinct traces.  One trace, the Yang-Mills measure,
is obtained by picking off the coefficient
of the empty skein.  The other four traces on $K_t(T^2)$ correspond to
each of the four $\mathbb{Z}_2$ homology
classes of the torus.
\end{abstract}

\maketitle


\section{Introduction}
Skein modules were introduced independently by Przytycki \cite{Przytycki1991}
and Turaev \cite{Turaev1988}
and have been an active topic of research since their introduction.
In particular, skein modules underlie quantum
invariants \cite{Lickorish1993, KauffmanLins} and are connected to the representation
theory of the fundamental group of the manifold \cite{Bullock1997, PrzytyckiSikora2000}.

The skein module is spanned by the equivalence classes of framed links
in the $3$-manifold.  The skein module
of the cylinder over a surface has a multiplication that
comes from laying one framed link on top of the other.  With this
multiplication, the skein module of the cylinder over a surface
is an algebra.

In this paper, we will consider the skein algebra of the torus specialized
at a complex number and describe the distinct traces on this
vector space.

\section{Preliminaries}

Let $M$ be an orientable 3-manifold.  A framed link
in $M$ is the embedding of disjoint annuli into $M$.  A framed
link is dipicted by drawing the core of each annulus.  One typically
uses the blackboard framing to produce the annulus from its core.
We will use $M = F \times I$ for a surface $F$.
In these cases, we will use the framing given by the
surface to produce the annulus from its core.

Equivalence of framed links in $M$ is up to regular isotopy.
That is, using only isotopy and Reidemeister's II and III moves.
A Reidemeister I move corresponds to a twist in the annulus and thus
such a move does not preserve the equivalence class of a framed link.

Let $\mathcal{L}(M)$ denote the equivalence class of framed links in $M$,
including the empty link, $\phi$.  Let $R = \mathbb{Z}[t,t^{-1}]$ be the ring
of Laurent polynomials.
Consider the free module
$R\mathcal{L}(M)$, with basis $\mathcal{L}(M)$.  Define $S(M)$ to
be the smallest subspace of $R\mathcal{L}(M)$ containing
all expressions of the form 
$\displaystyle{\lcr-t\zer-t^{-1}\ift}$
and 
$\bigcirc+t^2+t^{-2}$,
where the framed links in each expression are identical outside
the region pictured in the diagrams.  The
\emph{Kauffman bracket skein module}
$K(M)$ is the quotient $R\mathcal{L}(M)/S(M).$

Because $K(M)$ is defined using local relations on framed
links, two homeomorphic manifolds have
isomorphic skein modules.  Thus $K(M)$ is an invariant of the
$3$-manifold $M$.

\begin{figure}[ht]
  \begin{center}
    \leavevmode
    \epsfxsize = 10cm
    \epsfysize = 2cm
    \epsfbox{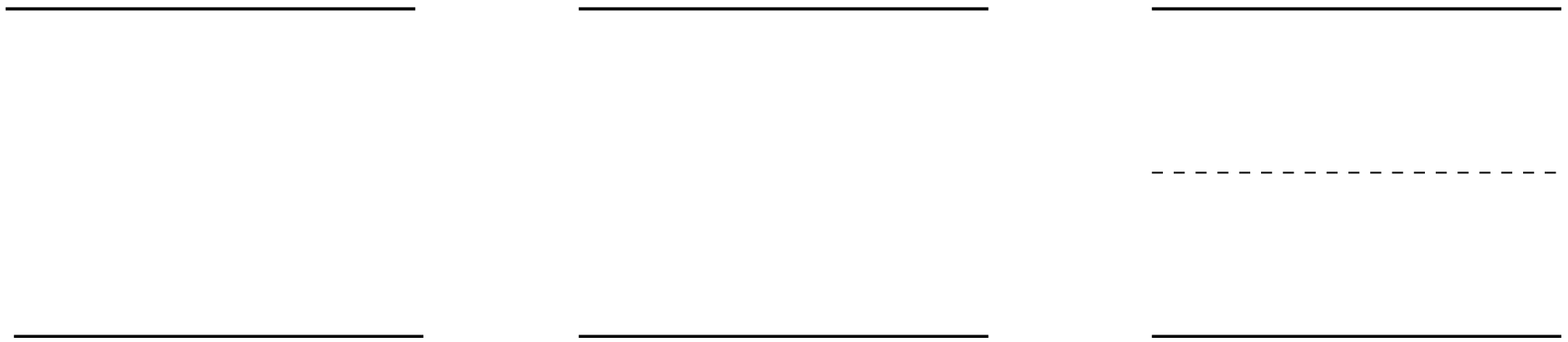}
    \put(-250,25){$\alpha$}
    \put(-200,25){$*$}
    \put(-145,25){$\beta$}
    \put(-93,25){$=$}
    \put(-45,35){$\alpha$}
    \put(-45,15){$\beta$}
    \caption{The product structure on $K(F \times I)$}
    \label{prod-struct}
  \end{center}
\end{figure}

Let $F$ be a compact, orientable surface and let
$I=[0,1]$ be the unit interval.  $K(F \times I)$
has an algebra structure that comes from laying one link
on top of the other.  Given skein elements $\alpha, \beta \in K(F \times I)$,
we can represent $\alpha, \beta$ with the links
$L_{\alpha}, L_{\beta} \in F \times I$.  Use isotopy to
move $L_{\alpha}$ to $F \times (\frac{1}{2},1]$ and
$L_{\beta}$ to $F \times [0,\frac{1}{2})$.  The
product $\alpha * \beta$ is the skein element
represented by $L_{\alpha} \cup L_{\beta}$.  A schematic of
this product structure is shown in Figure \ref{prod-struct}.

\vspace{0.5cm}

To simplify notation, and to emphasize that the algebra structure
is determined by $F$ rather than by $F \times I$, we will use
$K(F)$ to denote the skein module $K(F \times I)$ with the
algebra structure described above.  We will refer to $K(F)$ as
the \textit{skein algebra} of the surface $F$.

With the relations used to define $S(M)$, we can represent
any skein element as the linear combination of \emph{simple} diagrams, these are
diagrams with no crossings and no trivial components.  The
skein elements induced by such diagrams form a basis for $K(F)$.
This fact was observed by Bullock, Frohman, and
Kania-Bartoszy\'{n}ska in \cite{BullockFrohmanKania1999}.

In this paper we will make a number of simplifying assumptions.
Namely,  $F$ is the
standard 2-torus $T^2$ and 
$t$ is a complex number that is nonzero and not a root of unity.
Thus the polynomials in $R=\mathbb{Z}[t,t^{-1}]$ are evaluated at the
complex number $t$ and the skein
algebra $K(F)$ is \textit{specialized} at $t$ to form
$K_t(F)$, a vector space
over $\mathbb{C}$.  Throughout this paper, we will divide by
expressions of the form $(t^n - t^{-n})$.  Choosing $t$ to be
a complex number that is not a root of unity
allows this type of division without requiring
the use of rational functions.

Let $C_t(F) = [K_t(F), K_t(F)]$ be the vector space over
$\mathbb{C}$ with basis consisting of the commutators on $K_t(F)$.
A \emph{trace} on the algebra $K_t(F)$ is a linear functional
$\varphi : K_t(F) \to \mathbb{C}$ satisfying
$\varphi (\alpha * \beta) = \varphi (\beta * \alpha)$ for
all $\alpha, \beta \in K_t(F)$.  Since a trace is linear,
this condition can also be written
$\varphi (\alpha * \beta - \beta * \alpha) = 0$.
Then a trace $\varphi$ on $K_t(F)$ has $C_t(F) \subset {\rm ker}(\varphi)$.
So $\varphi$ descends to be a linear functional on the
quotient $K_t(F)/C_t(F)$.

\section{Examples of skein algebras}
\label{skein-alg-ex}

Before we explore the traces on $K_t(T^2)$ in detail, let's look a a few
examples of skein algebras.

Let $F$ be the 2-dimensional disk $D^2$.
Since every diagram in $D^2$ that has no crossings is
trivial, the only simple diagram in $D^2$ is the empty
skein $\phi$.  Thus the skein algebra $K(D^2)$ is one-dimensional
with basis $\{ \phi \}$.

Let $F$ be the annulus $A = S^1 \times [0,1]$.  The simple diagrams of $A$
consist of the empty skein along with any number of
parallel copies of the core of the annulus.  Denote $n$ parallel
copies of the core of the annulus by $z^n$ with $z^0 = \phi$.
Then a basis for $K(A)$ is $\{ z^0, z^1, z^2, \dots \}$ and
hence $K(A)$ is isomorphic to the algebra of polynomials in
$z$ with coefficients from $R = \mathbb{Z}[t, t^{-1}]$.

Now let $F$ be the torus $T^2$.  The collection of simple diagrams
has a more intricate structure.  We can have any $(p,q)$-curve with $p$ and $q$
relatively prime and also any number of parallel copies of a
given $(p,q)$-curve.  Thus a basis for $K(T^2)$ is
$\{ \phi, (p,q)^n ~\vert~ {\rm gcd}(p,q)=1, n \in \mathbb{N} \}$,
and we use the convention that $(p,q)^0 = \phi$.

Properties of the algebra $K(T^2)$ are explored
by Frohman and Gelca in \cite{FrohmanGelca2000}.  In particular,
they give a basis for $K(T^2)$ that behaves nicely under multiplication.
Let $(p,q)_T$ be the $(p,q)$-curve when ${\rm gcd}(p,q) = 1$ and let
\[
(p,q)_T = T_{ {\rm gcd}(p,q) } \left( \left(
\frac{p}{{\rm gcd}(p,q)}, \frac{q}{{\rm gcd}(p,q)}
\right) \right)
\]
when ${\rm gcd}(p,q) \neq 1$.  Here $T_n(x)$ is the $n$-th
Chebyshev polynomial defined recursively by $T_0(x) = 2$,
$T_1(x) = x$, and $T_{n+1} (x) = T_n(x) T_1(x) - T_{n-1}(x)$.
Using this basis for $K(T^2)$, we have the following
\textit{product-to-sum} formula.

\begin{theorem}(Frohman-Gelca)
\[
(p,q)_T * (r,s)_T = t^{|^{p~q}_{r~s}|}
(p+r,q+s)_T + t^{-|^{p~q}_{r~s}|}
(p-r, q-s)_T
\]
where ${|^{p~q}_{r~s}|}$ is the determinant.
\label{prod-to-sum}
\end{theorem}

\section{Traces on $K_t(T^2)$}
\label{hh0-t2}

Let $t$ be a fixed complex number, nonzero and not a root of unity.
Let $A = K_t(T^2)$ and let $\delta: A \ot A \to A$ be defined by
$\delta(a \ot b) = ab - ba$.  The image of $\delta$
in $A$ is the subalgebra of $A$ generated by these
commutators.  Denote this subalgebra by $C(A)$.
To understand the nature of the traces on $A$, we
will focus our attention on the commutator quotient $A / C(A)$ because
the space of traces is dual to this quotient.

As a vector space
$A/C(A)$ is spanned by the cosets of all
$(p,q)_T \in A$.  To narrow this spanning set to
a basis, we use the product-to-sum formula.

As cosets, $(x,y)_T + C(A) = (z,w)_T + C(A)$ if and only if
$(x,y)_T - (z,w)_T \in C(A)$.  
$(x,y)_T - (z,w)_T \in C(A)$ if and only
if $(x,y)_T - (z,w)_T$ is equal to some linear combination
of commutators.  The simplest case would be if
\[
(x,y)_T - (z,w)_T = \lambda \Big( (p,q)_T * (r,s)_T
- (r,s)_T * (p,q)_T \Big)~~~~~\lambda \in \mathbb{C}.
\]
We can use the product-to-sum formula to find the integers
$p$, $q$, $r$, $s$ when we are given $x$, $y$, $z$, $w$.

\begin{lemma}
Pick a complex number $t$ that is nonzero and not a root of unity and let $A = K_t(T^2)$
with basis $\{ (p,q)_T \}$ as described above.  Then
$(x,y)_T + C(A) = (z,w)_T + C(A)$ if $x + z$ is even, $y + w$ is
even, and ${|^{x~y}_{z~w}|} \neq 0$.
\label{lemma1}
\end{lemma}

\proof
Suppose we have integers $x$, $y$,
$z$, $w$ such that $x+z$ is even, $y + w$ is even and
${|^{x~y}_{z~w}|} \neq 0$.
Let $p = \frac{x+z}{2}$,
$q = \frac{y+w}{2}$, $r = \frac{x-z}{2}$, and $s = \frac{y-w}{2}$.
Since 
${|^{x~y}_{z~w}|} \neq 0$,
elementary matrix operations lead to
${|^{p~q}_{r~s}|} \neq 0$.
Let $\alpha =
{|^{p~q}_{r~s}|} \neq 0$.
Then 
\begin{eqnarray}
(p,q)_T * (r,s)_T & = & t^{\alpha} (x,y)_T + t^{-\alpha} (z,w)_T \nonumber \\
(r,s)_T * (p,q)_T & = & t^{-\alpha} (x,y)_T + t^{\alpha} (-z,-w)_T \nonumber
\end{eqnarray}
Since orientation doesn't matter in $A$,
$(-z,-w)_T = (z,w)_T$ and therefore
$$(p,q)_T * (r,s)_T - (r,s)_T * (p,q)_T =
(t^{\alpha} - t^{-\alpha})((x,y)_T - (z,w)_T).$$
Since $\alpha \neq 0$ and $t$ is not a root of unity, we can divide by
$(t^{\alpha} - t^{-\alpha})$ to get
\[
(x,y)_T - (z,w)_T = \frac{1}{t^{\alpha} - t^{-\alpha}}
\Big( (p,q)_T * (r,s)_T - (r,s)_T * (p,q)_T \Big) \in C(A).
\]
Thus $(x,y)_T + C(A) = (z,w)_T + C(A)$ as cosets in $A/C(A)$.
\qed

In other words, if $x$ and $z$ have the same parity, $y$ and $w$ have
the same parity, and $(x,y)$ and $(z,w)$ are linearly independent,
then $(x,y)_T$ is equivalent to $(z,w)_T$ in $A/C(A)$.
This fact allows us to reduce our basis for $A/C(A)$
somewhat, and it suggests that the parity of $(p,q)$
will determine the class of $(p,q)_T$ in $A/C(A)$.

\begin{theorem}
Pick a fixed complex number $t$, nonzero and not a root of unity.
Let $A = K_t(T^2)$ and let $C(A)$ be the subalgebra of $A$ generated
by commutators.  Then
$A/C(A)$ is a five dimensional
vector space over $\mathbb{C}$.
\label{five-dim}
\end{theorem}

\proof
Motivated by Lemma \ref{lemma1}, we define a map
$$\varphi \colon
A \to \mathbb{C} \left\{ \phi, ee, eo, oe, oo \right\}$$
by
$$(p,q)_T \mapsto \left\{
{\begin{array}{ccc}
\phi & {\rm if} & p=0,~q=0 \\
ee & {\rm if} & p~{\rm even},~q~{\rm even} \\
eo & {\rm if} & p~{\rm even},~q~{\rm odd} \\
oe & {\rm if} & p~{\rm odd},~q~{\rm even} \\
oo & {\rm if} & p~{\rm odd},~q~{\rm odd}
\end{array}}
\right.$$
then extend linearly.

We need to show that
${\rm ker}(\varphi) = C(A)$.
Recall that $C(A)$ is the vector space generated by
$(p,q)_T * (r,s)_T - (r,s)_T * (p,q)_T$
for $(p,q)_T, (r,s)_T \in A$.

Choose $(p,q)_T, (r,s)_T \in A$.
To show that $C(A) \subset {\rm ker}(\varphi)$, it suffices
to show that $c =
(p,q)_T * (r,s)_T - (r,s)_T * (p,q)_T
\in {\rm ker}(\varphi)$.
\begin{eqnarray}
c & = & (p,q)_T * (r,s)_T - (r,s)_T * (p,q)_T \nonumber \\
& = & \Big( t^{|^{p~q}_{r~s}|} - t^{- |^{p~q}_{r~s}|} \Big)
\Big( (p+r, q+s)_T - (p-r, q-s)_T \Big) \nonumber
\end{eqnarray}

If $(p+r,q+s)_T = (0,0)_T$, then $(p,q) = -(r,s)$.
Hence $(p,q)_T = (r,s)_T$.  Thus $c = 0 \in {\rm ker}(\varphi)$.
If $(p-r,q-s)_T = (0,0)_T$, then $(p,q) = (r,s)$.  Hence
$(p,q)_T = (r,s)_T$.  Thus $c = 0 \in {\rm ker}(\varphi)$.

If neither $(p+r,q+s)_T$ nor $(p-r,q-s)_T$ is $(0,0)_T$, then
since $p+r$ and $p-r$ have the same parity and $q+s$ and $q-s$
have the same parity, we have
$\varphi \left( (p+r, q+s)_T \right) = \varphi
\left( (p-r,q-s)_T \right)$ which implies that $c \in {\rm ker}(\varphi)$.
Hence $C(A) \subset {\rm ker}(\varphi)$.

Now we show that ${\rm ker}(\varphi) \subset C(A)$.
Take $k \in {\rm ker}(\varphi)$.  So $k \in A$
and $\varphi(k) = 0$.  Then
\begin{eqnarray}
k & = & \sum_{{\rm finite}} \lambda_{(p,q)} (p,q)_T \nonumber \\
& = & \lambda_{(0,0)} (0,0)_T + \sum_{ee} \lambda_{(p,q)}
(p,q)_T + \dots + \sum_{oo} \lambda_{(p,q)} (p,q)_T
\label{parity-sum}
\end{eqnarray}
Here we are breaking the sum into five parts, according to the
parity of $(p,q)$.  In each of these five parts, the coefficients
must sum to zero since $k \in {\rm ker}(\varphi)$.
As a model for the other cases, we work the case where
$$k = \sum_{ee} \lambda_{(p,q)} (p,q)_T$$ and so
$$\sum_{ee} \lambda_{(p,q)} = 0.$$

Now, $$\sum_{ee} \lambda_{(p,q)} (p,q)_T$$ is a finite sum
and $(p,q)_T$ are all of even-even parity, and $(p,q) \neq (0,0)$.
Choose integers
$r$ and $s$ such that $(r,s)_T$ is of even-even parity and
$(r,s)$ is linearly independent to each of the $(p,q)$ in
the sum for $k$.  That is, $\frac{s}{r}$ is a rational slope
that is different from the finite number of rational slopes
$\frac{q}{p}$.
Now using Lemma \ref{lemma1} each $(p,q)_T = (r,s)_T$ in
the quotient $A/C(A)$.
Hence
\begin{eqnarray}
k & = & \sum_{ee} \lambda_{(p,q)} (p,q)_T \nonumber \\
& = & \left( \sum_{ee} \lambda_{(p,q)} (r,s)_T \right) + {\rm commutators}
\nonumber \\
& = & \left( \sum_{ee} \lambda_{(p,q)} \right) (r,s)_T + {\rm commutators}
\nonumber \\
& = & {\rm commutators,~since} \sum_{ee} \lambda_{(p,q)} = 0. \nonumber
\end{eqnarray}
Thus $k \in C(A)$.

We could repeat this process for each of the $eo$, $oe$, and $oo$
sums given in Equation \ref{parity-sum}.  Thus for a general
$k \in \mathrm{ker}(\varphi)$, we have $k \in C(A)$.
Hence $\mathrm{ker}(\varphi) \subset C(A)$.
Now $\mathrm{ker}(\varphi) = C(A)$ and we have
\[
A/C(A) \cong \mathbb{C} \left\{ \phi, ee, eo, oe, oo \right\}.
\]
Thus $A/C(A)$ is a
five dimensional vector space over $\mathbb{C}$.
\qed

Recall that a trace is a linear functional defined on $A$ that is zero on $C(A)$.
The space of traces is dual to the quotient $A / C(A)$.
Thus Theorem \ref{five-dim} implies that there are five traces
on $A = K_t(T^2)$.
There is a trace for each $\mathbb{Z}_2$ homology class of
$T^2$, with one more trace for the empty skein.  Each trace picks off
the coefficients of the basis elements in its corresponding
class.  We could denote these traces by $\varphi_{\phi}$,
$\varphi_{ee}$, $\varphi_{eo}$,
$\varphi_{oe}$, $\varphi_{oo}$.  The trace $\varphi_{\phi}$ is
what Bullock, Frohman, and Kania-Bartoszy\'{n}ska call
the {\em Yang-Mills measure} in \cite{BullockFrohmanKania2003}.

\section{Further Investigation}

The Yang-Mills measure, $\varphi_{\phi}$, is defined on the skein algebra
of any closed surface $F$.  We have shown that there are four additional traces
when $F = T^2$.
Frohman and Kania-Bartoszy\'{n}ska have connected $\varphi_{\phi}$
to the $SU(2)$-characters of $\pi_1(F)$ in \cite{FrohmanKania2004}
and gave a state-sum formula for computing $\varphi_{\phi}$.

It is natural to ask if the other four traces, $\varphi_{ee}$, $\varphi_{eo}$,
$\varphi_{oe}$, and $\varphi_{oo}$, could be used in a similar way to study
the skein algebra of the torus or to study the skein modules of
manifolds with a torus boundary.

The commutator quotient $A/C(A)$ can also be seen as the zeroth Hochschild
homology of the skein algebra of the torus.  Understanding the structure
of the commutator quotient may have an impact on the Hochschild homology of
genus one Heegaard splittings as explored by the author in \cite{McLendon2006}.

\bibliography{torus}		
\bibliographystyle{plain}

\end{document}